\documentclass[fleqn,a4paper,11pt]{article}
\usepackage{amsmath,amsfonts,calrsfs,amssymb,color}
\usepackage{amsthm}

\newtheorem{Theorem}{Theorem}[section] 

\newtheorem{Proposition}[Theorem]{Proposition}
\newtheorem{Lemma}[Theorem]{Lemma}
\newtheorem{Corollary}[Theorem]{Corollary}
\newtheorem{Remark}[Theorem]{Remark}

\makeatletter
\@addtoreset{equation}{section}

\makeatother

\def\R{\mathbb R}
\def\N{\mathbb N}

\def\E{\mathbb E}
\def\P{\mathbb P}

\def\ds{\displaystyle}

 \def\T{\mathbb{T}}
 \def\eps{\varepsilon}

\title{The Schr\"odinger equation with spatial white noise potential}

\author{Arnaud Debussche\\
\normalsize IRMAR, ENS Rennes, UBL, CNRS
\and Hendrik Weber\\
\normalsize University of Warwick}
 
\date{}

\begin{document}

\maketitle

\begin{abstract}
We consider the linear and nonlinear Schr\"odinger equation with a spatial white noise as a potential in dimension $2$. We prove existence and uniqueness of solutions thanks to a change of unknown 
originally used in \cite{HL15} and conserved quantities.

\end{abstract}

\noindent {\bf 2010 Mathematics Subject Classification AMS}: \medskip

\noindent {\bf Key words}: Nonlinear Schr\"odinger, white noise.

\section{Introduction}

In this work we study a linear or nonlinear Schr\"odinger equation on a periodic domain with a random potential given by a spatial white noise in dimension $2$. This equation is important for various purposes. In the linear case, it is used to study Anderson localisation. It is a complex version of the famous PAM model. In the nonlinear case, it describes the evolution of 
nonlinear dispersive waves in a totally disorder medium (see for instance \cite{conti1}, \cite{conti2} and the references therein).  

If $u$ denotes the unknown, the equation is given by:
$$
\frac{d u }{dt}=\Delta u + \lambda |u|^{2} u + u \xi, \; x\in \T^2,\; t\ge 0,
$$
where $\T^2$ denotes the two dimensional torus, identified with $[0,2 \pi]^2$, and $\xi$ is a real-valued spatial white noise.  Of course, $\lambda=0$ for the linear equation. A positive $\lambda$ corresponds to the focusing case and $\lambda<0$ to the defocusing case. For simplicity, we take $\lambda=\pm 1$ in the nonlinear case. The qualitative properties of the solutions are completely different in these two cases. 

We are here interested in the question of existence and uniqueness of solutions. This is a preliminary but important step before studying other phenomena: solitary waves, blow up, Anderson 
localisation... 
The main difficulty is of course due to the presence of the rough potential. Recall that in dimension $2$, a white noise has a negative regularity which is strictly less than $-1$. Apparently, since the Schr\"odinger equation has very few smoothing properties and since this smoothing is very difficult to use, it seems hopeless to regularize it. 

However, this equation has many other properties. In particular, it is Hamiltonian and preserves the mass. Using a transformation due to \cite{HL15} in the context of PAM, we are able to use these invariants 
and construct solutions which have regularity strictly less than $1$. More precisely, we solve a transformed equation in almost $H^2(\T^2)$, the standard Sobolev space of functions with derivatives up to order $2$ in $L^2(\T^2)$. This is rather surprising since the regularity is comparable to what is obtained in the parabolic case, when strong smoothing properties are
available. 

As in the parabolic case, a renormalization is necessary and at the level of the original equation, the renormalized equation rewrites formally:
$$
\frac{d u }{dt}=\Delta u + \lambda |u|^{2} u + u (\xi-\infty), \; x\in \T^2,\; t\ge 0.
$$
The transformation $u\to e^{i\infty t}u$ transforms the original equation into the renormalized one. Therefore, the renormalization amounts to renormalize only the phase. A similar remark 
was made in \cite{bourgain} in a related but different context.  

For $\lambda<0$, we obtain global solution for any initial data satisfying some smoothness assumptions. For $\lambda>0$, as in the determinstic case, we need a smallness 
assumption on the initial data.

We could of course consider the equation with a more general nonlinearity: $|u|^{2\sigma}u$ with $\sigma\le 1$. For $\sigma<1$, no restriction on the size of the initial data is required for $\lambda>0$.  Another easy generalization is to consider a general bounded domain and Dirichlet boundary conditions, as long as they are sufficiently smooth and the properties of the 
Green function of the Laplace operator are sufficiently good so that Lemma \ref{l2.0} below holds. 

The study of the linear equation is closely related to the understanding of the Schr\"odinger operator with white noise potential. This is the subject of a recent very interesting article by Allez and Chouk (\cite{allez-chouk}) where the paracontrolled calculus is used to study the domain and spectrum of this operator. It is not clear how this can be used for the nonlinear equation.
\medskip

We use the classical $L^p=L^p(\T^2)$ spaces for $p\in [1,\infty]$, as well as the $L^2$ based Sobolev spaces $H^s=H^s(\T^2)$ for $s\in \R$ and the Besov spaces $B_{p,q}^s=B_{p,q}^s(\T^2)$, for $s\in\R, \, p,q\in [1,\infty]$.  These are defined
in terms of Fourier series and Littlewood-Paley theory (see \cite{BCD}).
Recall that 
$H^s=B^s_{2,2}$ and that, for $k\in\N, \; s\in(0,1)$, $B^{k+s}_{\infty,\infty}$ coincide with the H\"older space $C^{k,s}(\T^2)$.

Throughout the article, $c$ denotes a constant which may change from one line to the next. Also, we use a small parameter $0<\eps<e^{-1}$ and
$K_{\eps}$ is a random constant which can also change but such that $\E K_\eps^p $ is uniformly bounded in $\eps$ for all $p$. Similarly, for $0<\eps_1,\eps_2<e^{-1}$, the random constant $K_{\eps_1,\eps_2}$ may depend on $\eps_1,\eps_2$ but $\E K_{\eps_1,\eps_2}^p $ is uniformly bounded in $\eps_1, \eps_2$ for all $p$.

\section{Preliminaries}

We consider the following nonlinear Schr\"odinger equation in dimension $2$ on the torus, that is periodic boundary conditions are assumed, for the complex-valued unknown 
$u=u(x,t)$:
\begin{equation}\label{e1.1}
i\frac{d u }{dt}=\Delta u + \lambda |u|^{2} u + u \xi, \; x\in \T^2,\; t\ge 0.
\end{equation}
It is supplemented with  initial data 
$$
u(x,0)=u_0(x),\; x\in \T^2.
$$
We need some smoothness on the initial data, this will be made precise below. In the focusing case $\lambda>0$
 we need an extra assumption on the size of $ \|u_0\|_{L^2}$ which has be small enough (see \eqref{e4.3} below). 

The potential $\xi$ is random and is a real valued spatial white noise on $\T^2$. For simplicity, we assume that it has a zero spatial average. 
 The general case could be recovered by adding an additional Gaussian random potential which is constant in space. This would not change the analysis below.

Formally equation \eqref{e1.1} has two invariant quantities. Given a solution $u$ of \eqref{e1.1}, the mass:
$$
N(u(t))=\int_{\T^2} |u(x,t)|^2dx
$$
is constant in time as well as the energy:
$$
H(u(t))=\int_{\T^2} \frac12 |\nabla u(x,t)|^2 -\frac\lambda{4} |u(x,t)|^{4} -\frac12 \xi(x) |u(x,t)|^2 dx.
$$
This is formal because the noise $\xi$ is very rough. In dimension $1$, the noise has regularity $-1/2^-$ and belongs to $B_{\infty,\infty}^\alpha$ for any $\alpha<-1/2$, therefore the product $\xi |u|^2$ can be defined 
rigorously for $u\in H^1$ and this provides a bound in $H^1$. Existence and uniqueness through regularization of the noise and a compactness argument can then be obtained. 

In dimension $2$, the noise lives in any space with regularity $-1^{-}$, that is any regularity strictly less than $-1$, and the solution is not expected to be sufficiently smooth to compensate this. In fact, the product is almost well defined for $u\in H^1$ and we are in a situation similar to the two dimensional nonlinear heat equation with space time white noise. 
We expect that a renormalization is necessary. 
\medskip

Inspired by \cite{HL15}, we introduce: 
$$
Y=\Delta^{-1}\xi
$$
(note that this is well defined since we consider a zero average noise, we choose $Y$ also with a zero average) and $v=ue^{Y}$.  Then the equation for $v$ reads
\begin{equation}\label{e1.2}
i\frac{dv}{dt}=\Delta v -2\nabla v\cdot \nabla Y + v|\nabla Y|^2 + \lambda |v|^{2}v e^{-2 Y}.
\end{equation}
Now $Y$ has regularity $1^{-}$ and $\nabla Y$ is $0^{-}$. Thus this transformation has lowered the roughness of the most irregular term on the right hand side. At this point It is easier to see why we need a renormalization: the term $|\nabla Y|^2$ is not well defined since $\nabla Y$ is not a function. However, the roughness is mild here and it has been known for long that up to renormalization 
by a log divergent constant this square term can be defined in the second order Wiener chaos based on $\xi$.

Let us be more precise. Let $\rho_\eps=\eps^{-2}\rho(\frac{\cdot}\eps)$ be a compactly supported smooth mollifier and consider the smooth noise $\xi_\eps=\rho_\eps*\xi$.  We denote by $Y_\eps=\Delta^{-1}\xi_\eps$. Then it is proved in \cite{HL15} that for every $\kappa >0$, $\xi$ belongs almost surely to $B_{\infty,\infty}^{-1-\kappa}$ and,
as $\eps\to 0$, $\xi_\eps$ converges in probability to $\xi$ in $B_{\infty,\infty}^{-1-\kappa}$. 

Also, denoting by 
$$
C_\eps=\E\left(|\nabla Y_\eps|^2\right)
$$
the quantity $:|\nabla Y_\eps|^2:= |\nabla Y_\eps|^2-C_\eps$ converges in $L^p(\Omega;B_{\infty,\infty}^{-\kappa})$ for any $p\ge 1$, $\kappa>0$ to a random
variable $:|\nabla Y_\eps|^2:$ in the second Wiener chaos associated to $\xi$. It is easy to see that $C_\eps$ goes to $\infty$ as $|\ln \eps|$ as $\eps\to \infty$:
$$
\E\left(|\nabla Y_\eps|^2\right)\sim K_0   |\ln \eps |
$$
for some $K_0>0$. By stationarity  this quantity does not depend on $x\in\T^2$.

The precise result, whose proof can be found in \cite{HL15} in the more difficult case of the space variable in $\R^2$ (see the proofs of Lemma 1.1 and Proposition 1.3 in this work), is the following.
\begin{Lemma}\label{l2.0}
 $1\ge\kappa' > \kappa>0$ and any $p\ge 1$, there exist a constant $c$ independent of $\eps$ such that:
$$
\left[\E\left( \| Y_\eps -  Y\|_{B_{\infty,\infty}^{1-\kappa'}}^p\right)\right]^{\frac1p}\le c \eps^{\kappa -\frac2p}
$$
and 
$$
\left[\E\left( \|:|\nabla Y_\eps|^2: - :|\nabla Y|^2|\|_{B_{\infty,\infty}^{-\kappa'}}^p\right)\right]^{\frac1p}\le c \eps^{\kappa - \frac2p}.
$$
\end{Lemma}
\begin{Remark}
Using the monotonicity of stochastic $L^p$ norms in $p$,  one can drop the exponent $-\frac2p$ in the right hand side. We state the result in this way because 
this is the bound that one actually proves. Below, we use this bound with $\kappa-\frac2p=\frac{\kappa'}2$ (and $\kappa$ instead of $\kappa'$). 
\end{Remark}

Note that for $s<\tilde s$, $p,r\ge 1$, we have $B_{\infty,\infty}^{\tilde s}\subset B_{p,r}^s$.  Thus, bounds in the latter Besov spaces follow.
\medskip

Instead of solving equation \eqref{e1.2} for $v_\eps$, we consider:
\begin{equation}\label{e1.3}
i\frac{dv_\eps}{dt}=\Delta v_\eps -2\nabla v_\eps\cdot \nabla Y_\eps + v_\eps :|\nabla Y_\eps|^2: + \lambda |v_\eps|^{2}v_\eps e^{-2Y_\eps}
\end{equation}
and setting $u_\eps=v_\eps e^{-Y_\eps}$:
$$
i\frac{d u_\eps }{dt}=\Delta u_\eps + \lambda |u_\eps|^{2} u_\eps + (u_\eps-C_\eps)\xi_\eps, \; x\in \T^2,\; t\ge 0.
$$
Since $Y_\eps$ is smooth, it is classical to prove that these equations have a unique solution in $C([0,T];H^k)$ for an initial data in 
$H^k$, $k=1,\, 2$, provided the $L^2$ norm is small for $\lambda>0$ (see for instance \cite{cazenave}, Section 3.6). More details 
are given in Section \ref{s4}. 

The mass and energy are transformed into the two following quantities which are invariant under the dynamic for $v_\eps$:
$$
\tilde N_\eps(v_\eps(t))=\int_{\T^2} |v_\eps(x,t)|^2 e^{-2 Y_\eps(x)} dx
$$
and 
$$
\tilde H_\eps(v_\eps(t))= \int_{\T^2} \left(\frac12 |\nabla v_\eps (x,t)|^2+\frac12v_\eps^2 :|\nabla Y_\eps|^2:-\frac\lambda{4} |v_\eps(x,t)|^{4}e^{-2 Y_\eps(x)}\right)e^{-2 Y_\eps(x)}dx .
$$
Since the most irregular term $:|\nabla Y_\eps|^2:$ here is not as rough as $\xi$, this transformed energy is a much better quantity than the original one. It is possible to give a meaning to it for $\eps=0$ and 
use it to get bounds in  $H^1$.

Below, we use the following simple results.
\begin{Lemma}\label{l2.1bis}
For any $\kappa\in (0,1)$ and any $p\ge 1$, there exist a constant independent on $\eps$ such that:
$$
\E\left(\| e^{-2Y_\eps}-e^{-2Y}\|_{B^{1-\kappa}_ {\infty,\infty}}^p\right)\le c \eps^{-\frac2p+\kappa}.
$$
\end{Lemma}
{\bf Proof:} 
Since $B^{1-\kappa}_{\infty,\infty}$ is equal to the H\"older space $C^{1-\kappa}(\T^2)$ we have:
$$
\| e^{-2Y_\eps}-e^{-2Y}\|_{B^{1-\kappa}_ {\infty,\infty}}= \| e^{-2Y}(e^{-2(Y_\eps-Y)}-1)\|_{B^{1-\kappa}_ {\infty,\infty}}\le  \| e^{-2Y}\|_{B^{1-\kappa}_ {\infty,\infty}}\| e^{-2(Y_\eps-Y)}-1\|_{B^{1-\kappa}_ {\infty,\infty}}.
$$ 
Then we write:
$$
\begin{array}{l}
\ds \| e^{-2Y}\|_{B^{1-\kappa}_ {\infty,\infty}}\le 2\| e^{-2Y}\|_{L^\infty} \| Y\|_{B^{1-\kappa}_ {\infty,\infty}}, \\
\ds \| e^{-2(Y_\eps-Y)}-1\|_{B^{1-\kappa}_ {\infty,\infty}}\le 2\| e^{-2Y}\|_{L^\infty}\| e^{-2Y_\eps}\|_{L^\infty} \|Y_\eps-Y\|_{B^{1-\kappa}_ {\infty,\infty}}
\end{array}
$$
The result follows by H\"older inequality, Lemma \ref{l2.0} and Gaussianity  to bound exponential moments of $Y$ and $Y_\eps$.
\hfill$\square$
\begin{Lemma}\label{l2.1}
There exists a contant $c$ independent of $\eps$ such that:
$$
\E\left(\|\nabla Y_\eps\|_{L^4}^4\right)\le c |\ln \eps|^2
$$
and 
$$
\E\left(\|:|\nabla Y_\eps|^2:\|_{L^4}^4\right)\le c (|\ln \eps|)^4.
$$
\end{Lemma}
{\bf Proof:} 
It suffices to write:
$$
\E\left(\int_{\T^2} |\nabla Y_\eps(x) |^4 dx \right)= \int_{\T^2}\E\left( |\nabla Y_\eps(x) |^4 \right) dx= 12\pi^2 C_\eps^2.
$$
Similarly:
$$
\E\left(\int_{\T^2} |:|\nabla Y_\eps(x)|^2: |^4 dx \right)=\E\left(\int_{\T^2} (|\nabla Y_\eps(x)|^2 -C_\eps)^4 dx \right)= 51\pi C_\eps^4.
$$
\hfill$\square$

\section{The linear case}\label{s3}

In this  section, we start with the linear case: $\lambda=0$. Then the equation for $v_\eps$ reads
\begin{equation}
\label{e3.1}
i\frac{dv_\eps}{dt}=\Delta v_\eps -2\nabla v_\eps\cdot \nabla Y_\eps + v_\eps :|\nabla Y_\eps|^2: .
\end{equation}
There exists a unique solution in $C([0,T];H^2)$ if $v_\eps(0)\in H^2$.  We 	take the initial data 
$$
v_\eps(0)=v_0=u_0e^{Y}
$$ 
and assume below that it belongs to $H^2$. 

The mass and energy of a solution are now:
$$
\tilde N_\eps(v_\eps(t))=\int_{\T^2} |v_\eps(x,t)|^2 e^{-2 Y_\eps(x)} dx
$$
and 
$$
\tilde H_\eps(v_\eps(t))= \int_{\T^2} \left(\frac12 |\nabla v_\eps (x,t)|^2+\frac12v_\eps^2 :|\nabla Y_\eps|^2:\right)e^{-2 Y_\eps(x)}dx .
$$
They are constant in time under the evolution \eqref{e3.1}.

Since $Y_\eps$ converges in $B_{\infty,\infty}^{1-\kappa}$ for any $\kappa>0$ as $\eps$ tends to zero, we see that the mass gives a uniform bound in 
$L^2$ on $v_\eps$. 
More precisely:
\begin{equation}\label{e3.2}
\|v_\eps(t)\|_{L^2}^2\le \|e^{2Y_\eps}\|_{L^\infty}\|e^{-2Y_\eps}\|_{L^\infty}  \|v_0\|_{L^2}^2 = K_\eps \|v_0\|_{L^2}^2 
\end{equation}
with
\begin{equation}\label{e3.3}
K_{\eps}=\|e^{2Y_\eps}\|_{L^\infty}\|e^{-2Y_\eps}\|_{L^\infty} .
\end{equation}

The energy enables us to get a bound on the gradient.
\begin{Proposition}\label{p3.1}
Let $\kappa\in (0,1/2)$, there exists a random constant $K_{\eps}$ bounded in $L^p(\Omega)$ with respect to $\eps$ for any $p\ge 1$  such that if $v_0\in H^1$:
$$
\int_{\T^2}  |\nabla v_\eps (x,t)|^2 dx \le K_{\eps}  \left( \tilde H_\eps(v_0) +  \|v_0\|_{L^2}^2
\right).$$
\end{Proposition}
{\bf Proof:} 
Since $B^{-\kappa}_{\infty,2}$ is in duality with $B^\kappa_{1,2}$ we deduce by the standard multiplication rule in Besov spaces (see e.g. \cite{BCD}, Section 2.8.1)
$$
\left|\int_{\T^2} v_\eps^2 :|\nabla Y_\eps|^2: dx\right|\le  \|v_\eps^2\|_{B^\kappa_{1,2}}\| :|\nabla Y_\eps|^2: \|_{B^{-\kappa}_{\infty,2}}
\le K_{\eps} \|v_\eps^2\|_{B^\kappa_{1,2}}\le  K_{\eps} \|v_\eps\|_{B^{\kappa}_{2,2}}^2.
$$
Then we note that $\|v_\eps\|^2_{B^{\kappa}_{2,2}}=\|v_\eps\|_{H^{\kappa}}^2$ so that by interpolation
$$
\left|\int_{\T^2} v_\eps^2 :|\nabla Y_\eps|^2: dx\right|\le  K_{\eps} \|v_\eps\|_{L^2}^{2(1-\kappa)}\|v_\eps\|_{H^1}^{2\kappa}.
$$
It follows
$$
\begin{array}{ll}
\ds \int_{\T^2}  |\nabla v_\eps (x,t)|^2 dx &\ds \le K_{\eps} \tilde H(v_\eps(t)) + K_{\eps} \|v_\eps(t)\|_{L^2}^{2(1-\kappa)}\|v_\eps(t)\|_{H^1}^{2\kappa}\\
&\ds \le  K_{\eps}  \tilde H_\eps(v_0) +  K_{\eps}   \|v_\eps(t)\|_{L^2}^2 + K_{\eps}   \|v_\eps(t)\|_{L^2}^{2(1-\kappa)}\|\nabla v_\eps\|_{L^2}^{2\kappa}\\
&\ds \le K_{\eps}  \tilde H_\eps(v_0) + K_{\eps}  \|v_0\|_{L^2}^2+\frac 12\|\nabla v_\eps\|_{L^2}^2
\end{array}
$$
and hence, by absorbing the last term in the left hand side,
$$
\int_{\T^2}  |\nabla v_\eps (x,t)|^2 dx \le K_{\eps} \left(  \tilde H_\eps(v_0) +   \|v_0\|_{L^2}^2\right).
$$
\hfill$\square$

Since $ \tilde H_\eps(v_0)$ is  bounded for $v_0\in H^1$, we obtain a (random) bound on $v_\eps$ in $H^1$ using similar arguments as above. 
\begin{Corollary}\label{c3.2}
There exists a random constant $K_{\eps}$ bounded in $L^p(\Omega)$ with respect to $\eps$ for any $p\ge 1$  such for any $v_0\in H^1$:
$$
 \|v_\eps(t)\|_{H^1} \le K_{\eps}     \|v_0\|_{H^1}, \; t\ge 0.
$$
\end{Corollary}

Unfortunately, this regularity is not sufficient to  control the product $\nabla v_\eps \cdot \nabla Y_\eps$ on the right hand side of \eqref{e3.1}. 

The next observation is that, if $v_0$ is smooth enough, $v_\eps$ is time differentiable and setting $w_\eps=\frac{d v_\eps}{dt}$ it is a solution of:
\begin{equation}
\label{e3.2bis}
i\frac{dw_\eps}{dt}=\Delta w_\eps -2\nabla w_\eps\cdot \nabla Y_\eps + w_\eps :|\nabla Y_\eps|^2: . 
\end{equation}
Since $w_\eps$ satisfies the same equation as $v_\eps$, it has the same invariant quantities. We use in particular the mass:
$$
\tilde N(w_\eps(t))=\tilde N(w_\eps(0)).
$$
Hence:
$$
\|w_\eps(t)\|_{L^2}\le K_{\eps} \tilde N(w_\eps(0))\le K_\eps \|w_\eps(0)\|_{L^2}.
$$
This is still true under the assumption that $w_\eps(0)\in L^2$, which is equivalent to $v_0\in H^2$. 
\begin{Proposition}\label{p3.2}
There exist a random constant $K_{\eps}$ bounded in $L^p(\Omega)$ with respect to $\eps$ for any $p\ge 1$ such that if $v_0\in H^2$:
$$
\|v_\eps\|_{H^2}\le cK_{\eps}\left(\|v_0\|_{H^2}+ \|v_0\|_{L^2} |\ln \eps|^2\right).
$$
\end{Proposition}
{\bf Proof:}
From \eqref{e3.1}, we have: 
$$
w_\eps(0)= -i(\Delta v_0  -2\nabla v_0\cdot \nabla Y_\eps + v_0 :|\nabla Y_\eps|^2:),
$$
so that, thanks to the embedding $H^{1/2}\subset L^4$, 
$$
\|w_\eps(0)\|_{L^2} \le c\left(\|v_0\|_{H^2}+ \|v_0\|_{H^{3/2}}\|\nabla Y_\eps\|_{L^4}+
\|v_0 \|_{H^{1/2}}\|:|\nabla Y_\eps|^2:\|_{L^4}\right).
$$
By interpolation we deduce:
\begin{equation}\label{e3.5bis}
\begin{array}{ll}
\|w_\eps(0)\|_{L^2}& \le c\left(\|v_0\|_{H^2}+ \|v_0\|_{H^{2}}^{3/4}\|v_0\|_{L^2}^{1/4}\|\nabla Y_\eps\|_{L^4}+
\|v_0\|_{H^{2}}^{1/4}\|v_0\|_{L^2}^{3/4}\|:|\nabla Y_\eps|^2:\|_{L^4}\right)\\
\\
&\ds \le c\left(\|v_0\|_{H^2}+ \|v_0\|_{L^2} \|\nabla Y_\eps\|_{L^4}^4+
\|v_0\|_{L^2}\|:|\nabla Y_\eps|^2:\|_{L^4}^{4/3}\right)\\
\\
&\ds \le c\left(\|v_0\|_{H^2}+ K_{\eps}\|v_0\|_{L^2} |\ln \eps|^2\right) ,
\end{array}
\end{equation}
where: 
$$
K_{\eps}= \|\nabla Y_\eps\|_{L^4}^4 |\ln \eps|^{-2}+\|:|\nabla Y_\eps|^2:\|_{L^4}^{4/3} |\ln \eps|^{-2}.
$$
By Lemma \ref{l2.1} and gaussianity, we know that the moments  of this random variable are bounded  with respect to $\eps$.

It follows
$$
\|w_\eps(t)\|_{L^2}\le K_{\eps}  \left(\|v_0\|_{H^2}+ \|v_0\|_{L^2} |\ln \eps|^2\right) .
$$
This in turn allows us to control $\|v_\eps\|_{H^2}$. Indeed, from \eqref{e3.1},
$$
\begin{array}{ll}
\|\Delta v_\eps\|_{L^2}&\le \|w_\eps(t)\|_{L^2}+2\|\nabla v_\eps\cdot \nabla Y_\eps\|_{L^2} + \|v_\eps :|\nabla Y_\eps|^2:\|_{L^2}
\end{array}
$$
and by similar arguments as above
$$
\|\Delta v_\eps\|_{L^2}\le \|w_\eps(t)\|_{L^2}+\frac12\|\Delta v_\eps\|_{L^2} +c K_{\eps}\|v_\eps(t)\|_{L^2} |\ln \eps|^2
$$
and finally
$$
\|\Delta v_\eps\|_{L^2}\le  K_{\eps}\left(\|v_0\|_{H^2}+ \|v_0\|_{L^2} |\ln \eps|^2\right).
$$
The result follows thanks to \eqref{e3.2}.
\hfill$\square$

This bound does not seem to be very useful since it explodes as $\eps\to 0$. To use it, we consider the difference of two solutions.
\begin{Proposition}\label{p3.3}
Let $\eps_2>\eps_1>0$ then for $\kappa\in (0,1]$, $p\ge1$ there exist a random constant $K_{\eps_1,\eps_2}$ bounded in $L^p(\Omega)$ with respect to $\eps_1,\eps_2$ for any $p\ge 1$ such that if $v_0\in H^2$
$$
\sup_{t\in [0,T]} \|v_{\eps_1}(t)-v_{\eps_2}(t)\|_{L^2}^2\le K_{\eps_1,\eps_2} \eps_2^{\kappa/2}|\ln \eps_2|^{1+2\kappa}\|v_0\|_{H^2}^2
$$
\end{Proposition}
{\bf Proof:}
We set $r=v_{\eps_1}-v_{\eps_2}$ and write:
$$
i\frac{dr}{dt}=\Delta r-2r\cdot \nabla Y_{\eps_1} + r :|\nabla Y_{\eps_1}|^2:
-2\nabla v_{\eps_2}\cdot \nabla (Y_{\eps_1}-Y_{\eps_2}) + v_{\eps_2} (:|\nabla Y_{\eps_1}|^2: -:|\nabla Y_{\eps_2}|^2:).
$$
By standard computations, we deduce:
$$
\begin{array}{l}
\ds \frac12\frac{d}{dt}\int_{\T^2}|r(x,t)|^2 e^{-2Y_{\eps_1}(x)}dx \\
\ds = Im \int_{\T^2} \left(-2\nabla v_{\eps_2}\cdot \nabla (Y_{\eps_1}-Y_{\eps_2}) + v_{\eps_2} (:|\nabla Y_{\eps_1}|^2: -:|\nabla Y_{\eps_1}|^2:)\right)\bar r e^{-2Y_{\eps_1}(x)}dx\\
\ds \le 2\|\nabla v_{\eps_2}\bar r e^{-2Y_{\eps_1}}\|_{ B_{1,2}^\kappa}\|\nabla (Y_{\eps_1}-Y_{\eps_2}) \|_{ B_{\infty,2}^{-\kappa}}
+\|v_{\eps_2}\bar r e^{-2Y_{\eps_1}}\|_{ B_{1,2}^\kappa}\|:|\nabla Y_{\eps_1}|^2: -:|\nabla Y_{\eps_2}|^2: \|_{ B_{\infty,2}^{-\kappa}}
\end{array}
$$
the first term of the right hand side is bounded thanks to interpolation and paraproduct inequalities (see \cite{BCD}) and we have
$$
\begin{array}{ll}
\|\nabla v_{\eps_2}\bar r e^{-2Y_{\eps_1}}\|_{ B_{1,2}^\kappa}&\ds \le c \|v_{\eps_2}\|_{H^{1+\kappa}}\left(\|v_{\eps_1}\|_{H^{\kappa}}+\|v_{\eps_2}\|_{H^{\kappa}}\right)\|e^{-2Y_{\eps_1}}\|_{B^{\kappa}_{\infty,2}}.\\
\end{array}
$$
Then, by Proposition \ref{p3.2} and interpolation: 
$$
\|v_{\eps_2}\|_{H^{1+\kappa}}\le c \|v_{\eps_2}\|_{H^2}^{(1+\kappa)/2}\|v_{\eps_2}\|_{L^2}^{(1-\kappa)/2}
\le K_{\eps_2} (\|v_0\|_{H^2}+\|v_0\|_{L^2}|\ln \epsilon_2|^2)^{(1+\kappa)/2} \|v_{0}\|_{L^2}^{(1-\kappa)/2}
$$
and, by Corollary \ref{c3.2}, for $i=1,2$ 
$$
\|v_{\eps_i}\|_{H^{\kappa}}\le \|v_{\eps_i}\|_{H^{1}}^{\kappa} \|v_{\eps_i}\|_{L^2}^{1-\kappa}\le K_{\eps_i}\|v_0\|_{H^1}^{\kappa}\|v_0\|_{L^2}^{1-\kappa}
\le K_{\eps_i}\|v_0\|_{H^2}^{\kappa/2}\|v_0\|_{L^2}^{1-\kappa/2}.
$$
It follows
$$
\|\nabla v_{\eps_2}\bar r e^{-2Y_{\eps_1}}\|_{ B_{1,2}^\kappa}\le K_{\epsilon_2}\|e^{-2Y_{\eps_1}}\|_{B^{\kappa}_{\infty,2}}\|v_0\|_{L^2}^{\frac32-\kappa}(\|v_0\|_{H^2}+\|v_0\|_{L^2}|\ln \epsilon_2|^2)^{\frac12+\kappa}
$$

 The second term is bounded by the same quantity and we deduce:
$$
\begin{array}{l}
\ds \frac12\frac{d}{dt}\int_{\T^2}|r(x,t)|^2 e^{-2Y_{\eps_1}(x)}dx \\
\\
\ds \le K_{\epsilon_2}\|e^{-2Y_{\eps_1}}\|_{B^{\kappa}_{\infty,2}}\|v_0\|_{L^2}^{\frac32-\kappa}(\|v_0\|_{H^2}+\|v_0\|_{L^2}|\ln \epsilon_2|^2)^{\frac12+\kappa}\\
\ds\hspace{1.5cm} \times
\left(\|\nabla (Y_{\eps_1}-Y_{\eps_2}) \|_{ B_{\infty,2}^{-\kappa}}+\|:|\nabla Y_{\eps_1}|^2: -:|\nabla Y_{\eps_2}|^2: \|_{ B_{\infty,2}^{-\kappa}}\right).
\end{array}
$$
The result follows thanks to integration in time and Lemma \ref{l2.0}. 
\hfill$\square$

By interpolation, we deduce from Proposition \ref{p3.2} and \ref{p3.3} the following result.
\begin{Corollary}\label{c3.4}
Let $\eps_2>\eps_1>0$ then for $\kappa\in (0,1]$, $\gamma\in [0,2)$, $p\ge1$ there exist a random constant $K_{\eps_1,\eps_2}$ bounded in $L^p(\Omega)$ with respect to $\eps_1,\eps_2$ for any $p\ge 1$ such that if $v_0\in H^2$:
$$
\sup_{t\in [0,T]}\|v_{\eps_1}(t)-v_{\eps_2}(t)\|_{H^\gamma}^2\le K_{\eps_1,\eps_2} \eps_2^{\frac\kappa2(1-\frac\gamma2)}|\ln \eps_1|^{4}\|v_0\|_{H^2}^2
$$
\end{Corollary}

We are now ready to state and prove the main result of this section.
\begin{Theorem}\label{t3.1}
Assume that $v_0=u_0e^{Y}\in  L^p(\Omega;H^2)$ for some $p>1$. For any $T\ge 0$, $q<p$, $\gamma\in (1,2)$, when $\eps\to0$, the solution $v_\epsilon$ of \eqref{e3.1} satisfying $v_\eps(0)=v_0$ converges in $L^q(\Omega;C([0,T];H^\gamma))$ to $v$ which is the unique solution to 
\begin{equation}
\label{e3.6}
i\frac{dv}{dt}=\Delta v -2\nabla v\cdot \nabla Y + v :|\nabla Y|^2: 
\end{equation}
in $C([0,T];H^\gamma)$ such that $v(0)=v_0$.
\end{Theorem} 
{\bf Proof:} Pathwise uniqueness is clear. Indeed if $v\in C([0,T];H^\gamma)$ is a solution of \eqref{e3.6}, it is not difficult to use a regularization argument 
to justify
$$
\frac{d}{dt}\int_{\T^2} |v(x,t)|^2e^{-2Y(x)}dx= 0.
$$
Now let $\eps_k=2^{-k}$,  Corollary \ref{c3.4} implies that  $(v_{\eps_{k}})$ is Cauchy in $L^q(\Omega;C([0,T];H^\gamma))$. 
It is not difficult to prove that the limit $v$ is a solution of \eqref{e3.6} and 
$$
 \E\left(\sup_{t\in [0,T]}\|v_{\eps_k}(t)\|_{H^\gamma}^{q}\right)\le c \E\left(\|v_0\|_{H^2}^p\right)^{q/p},
$$ 
$$
 \E\left(\sup_{t\in [0,T]}\|v(t)\|_{H^\gamma}^{q}\right)\le c \E\left(\|v_0\|_{H^2}^p\right)^{q/p}.
$$ 
Then, by interpolation with $\gamma<\tilde \gamma  <2$: 
$$
\begin{array}{ll}
&\ds \sup_{t\in [0,T]}\|v_{\eps}(t)-v_{\eps_k}(t)\|_{H^\gamma}^{2}\ds\\
&\ds \le c \sup_{t\in [0,T]}\|v_{\eps}(t)-v_{\eps_k}(t)\|_{L^2}^{2(1-\gamma/{\tilde\gamma})} \sup_{t\in [0,T]}\|v_{\eps}(t)-v_{\eps_k}(t)\|_{H^{\tilde\gamma}}^{2\gamma/{\tilde\gamma}}\\
&\ds \le c \sup_{t\in [0,T]}\|v_{\eps}(t)-v_{\eps_k}(t)\|_{L^2}^{2(1-\gamma/{\tilde\gamma})} \sup_{t\in [0,T]}\left(\|v_{\eps}(t)\|_{H^2}+\|v_{\eps_k}(t)\|_{H^{\tilde\gamma}}\right)^{2\gamma/{\tilde\gamma}}.
\end{array}
$$
By Proposition \ref{p3.3}, Proposition \ref{p3.2} and the above inequality, we deduce:
$$
\E\left(\sup_{t\in [0,T]}\|v_{\eps}(t)-v_{\eps_k}(t)\|_{H^\gamma}^{q}\right)\le c \eps^{\frac\kappa2(2-\gamma/{\tilde\gamma}))}|\ln \eps|^{4}\E\left(\|v_0\|_{H^2}^p\right)^{q/p}.
 $$
Letting $k\to \infty$, we deduce that the whole family $(v_\eps)_{\eps>0}$ converges to $v$ in $L^q(\Omega;C([0,T];H^\gamma))$.
\hfill$\square$

\section{The nonlinear equation}\label{s4}

We now study the nonlinear equation \eqref{e1.2} and consider its approximation \eqref{e1.3} with initial condition
$$
v_\eps(0)=v_0=u_0e^{Y}
$$ 
and assume below that it belongs to $H^2$. 

 Again the mass gives a uniform bound in 
$L^2$ on $v_\eps$:
\begin{equation}\label{e4.1}
\|v_\eps(t)\|_{L^2}^2\le K_{\eps}  \|v_0\|_{L^2}^2.
\end{equation}

The estimate on the $H^1$ norm using the energy is similar to the linear case. Recall that the energy is given by:
$$
\tilde H_\eps(v_\eps(t))= \int_{\T^2} \left(\frac12 |\nabla v_\eps (x,t)|^2+\frac12v_\eps^2 :|\nabla Y_\eps|^2:-\frac\lambda{4} |v_\eps(x,t)|^{4}e^{-2 Y_\eps(x)}\right)e^{-2 Y_\eps(x)}dx 
$$
and it can be checked that for all $t\geq 0$ we have $\tilde H_\eps(v_\eps(t))=\tilde H_\eps(v_0)$.
\begin{Proposition}\label{p4.1}
There exists a constant $K_{\eps}$ bounded in $L^p(\Omega)$ for any $p\ge 1$  such that if $v_0\in H^1$ and 
\begin{equation}\label{e4.2}
\|e^{-2Y_\eps}\|_{L^\infty}^3\|e^{2Y_\eps}\|_{L^\infty} \|v_0\|_{L^2}\le 1 \mbox{ if } \lambda=1
\end{equation}
then
$$
\|v_\eps (t)\|_{H^1}^2 \le K_{\eps}  \left( \|v_0\|_{H^1}^2+ \|v_0\|_{L^2}^2\|v_0\|_{H^1}^2\right).
$$
Moreover, the sequence $(K_k)=(K_{2^{-k}})$ is bounded almost surely.
\end{Proposition}
\begin{Remark}
It is classical that the equation for $u_\eps$ is locally well posed in $H^1$, see \cite{cazenave} Section 3.6. The bound obtained here gives a global bound in $H^1$ for $u_\eps$ so that 
global existence and uniqueness hold for $u_\eps$ and thus for $v_\eps$.
\end{Remark}
{\bf Proof:} 
We proceed as in the proof of Proposition \ref{p3.1}. We first have
$$
\left|\int_{\T^2} v_\eps^2 :|\nabla Y_\eps|^2: dx\right|\le  K_{\eps} \|v_\eps\|_{L^2}^{2(1-\kappa)}\|v_\eps\|_{H^1}^{2\kappa}
$$
and 
$$
\begin{array}{ll}
&\ds \int_{\T^2}  |\nabla v_\eps (x,t)|^2 dx \\
&\ds \le K_{\eps} \tilde H(v_\eps(t)) + K_{\eps} \|v_\eps(t)\|_{L^2}^{2(1-\kappa)}\|v_\eps(t)\|_{H^1}^{2\kappa}
+\frac\lambda{4} \int_{\T^2} |v_\eps(x,t)|^{4}e^{-2 Y_\eps(x)}e^{-2 Y_\eps(x)}dx \\
&\ds \le K_{\eps}  \tilde H_\eps(v_0) + K_{\eps}  \|v_0\|_{L^2}^2+\frac 12\|\nabla v_\eps\|_{L^2}^2
+\frac\lambda{4} \int_{\T^2} |v_\eps(x,t)|^{4}e^{-2 Y_\eps(x)}e^{-2 Y_\eps(x)}dx .
\end{array}
$$

For $\lambda=-1$, the result follows after dropping the last term and using
$$
\int_{\T^2} |v_0(x)|^4 e^{-4Y_\eps(x)}dx \le K_\eps \|v_0\|_{L^4}^4 \le K_\eps \|v_0\|_{H^{1/2}}^4\le K_\eps  \|v_0\|_{L^2}^{2}\|v_0\|_{H^1}^{2}
$$
thanks to the Sobolev embedding $H^{1/2}\subset L^4$ and interpolation.

For $\lambda=1$,  Gagliardo-Nirenberg inequality (see for instance \cite{brezis-gallouet} for a simple proof with the constant $1/2$ used below): 
$$
\begin{array}{ll}
\ds \int_{\T^2} |v_\eps(x,t)|^4 e^{-4Y_\eps(x)}dx &\ds \le \|e^{ - 2Y_\eps}\|_{L^\infty}^2 \int_{\T^2} |v_\eps(x,t)|^4 dx\\
&\ds  \le \frac12 \|e^{-2Y_\eps}\|_{L^\infty}^2 \int_{\T^2} |v_\eps(x,t)|^2 dx \int_{\T^2} |\nabla v_\eps(x,t)|^2 dx\\
&\ds \le \frac12 \|e^{-2Y_\eps}\|_{L^\infty}^3 \|e^{2Y_\eps}\|_{L^\infty} \|v_0\|_{L^2}^2 \int_{\T^2} |\nabla v_\eps(x,t)|^2 dx,
\end{array}
$$
where we have made use of $\|v_\eps(t)\|_{L^2}^2\le \|e^{2Y_\eps}\|_{L^\infty}\|e^{-2Y_\eps}\|_{L^\infty}  \|v_0\|_{L^2}^2 $ according to \eqref{e3.2}. The result follows easily under assumption \eqref{e4.2}. 
 
 The constant $K_\eps$ is  a polynomial in $\| :|\nabla Y_\eps|^2: \|_{B^{-\kappa}_{\infty,2}}$, $\|e^{-2Y_\eps}\|_{L^\infty}$ and $\|e^{2Y_\eps}\|_{L^\infty}$.
  By Lemma \ref{l2.0} and Lemma \ref{l2.1bis} and Borel-Cantelli, we know that $:|\nabla Y_{2^{-k}}|^2:$ and $Y_{2^{-k}}$ converge almost surely in 
 $B^{-\kappa}_{\infty,\infty}$ and $B^{1-\kappa}_{\infty,\infty}$ so that $K_k$ is indeed bounded almost surely. 
 \hfill$\square$

We now proceed with the $H^2$ bound.
\begin{Proposition}\label{p4.2}
There exist a random constants $K_{\eps}$ bounded in $L^p(\Omega)$ with respect to $\eps$ for any $p\ge 1$ such that if $v_0=u_0e^{-Y}\in H^2$ 
and \eqref{e4.2} holds:
$$\begin{array}{l}
\|v_\eps(t)\|_{H^2}\\ \le 
cK_{\eps}\left(1+\|v_0\|_{H^2}+ \|v_0\|_{L^2} |\ln \eps|^4+\|v_0\|_{H^1}^3+\|v_0\|_{L^2}^3\|v_0\|_{H^1}^3\right)^{\exp(K_{\eps}\left(\|v_0\|_{H^1}^2+\|v_0\|_{L^2}^2\|v_0\|_{H^1}^2\right)t)}.
\end{array}
$$
Moreover $K_k=K_{2^{-k}}$ is bounded almost surely.
\end{Proposition}
\begin{Remark}
We know that if $u_\eps(0)\in H^2$, the solution $u_\eps$ lives in $H^2$ (see \cite{brezis-gallouet}). Since $u_\eps(0)=v_0e^{-Y_\eps}$, we can apply this result. 
\end{Remark}
{\bf Proof:} As in Section 3, we set $w_\eps=\frac{dv_\eps}{dt}$ which now satisfies:
$$
i\frac{dw_\eps}{dt}=\Delta w_\eps-2\nabla w_\eps\cdot \nabla Y_\eps +w_\eps :|\nabla Y|^2: + \lambda \left(|v_\eps|^2 w_\eps + 2 Re(v_\eps \bar w_\eps) v_\eps \right) e^{-2Y_\eps}.
$$

From \eqref{e1.3}, we have: 
$$
w_\eps(0)= -i(\Delta v_0  -2\nabla v_0\cdot \nabla Y_\eps + v_0 :|\nabla Y_\eps|^2:)+\lambda |v_0|^2v_0e^{-2Y_\eps},
$$
and as in Proposition \ref{p3.2} and using the embedding $H^{1}\subset L^6$: 
$$
\|w_\eps(0)\|_{L^2}\ds \le c \|v_0\|_{H^2}+ \|v_0\|_{L^2}\left(\|\nabla Y_\eps\|_{L^4}^4+\|:|\nabla Y_\eps|^2:\|_{L^4}^{4/3}   \right)+\|e^{-2Y_\eps}\|_{L^\infty}\|v_0\|_{H^1}^3.
$$
By Lemma \ref{l2.1}:
$$
\P(\|\nabla Y_\eps\|_{L^4}^4+\|:|\nabla Y_\eps|^2:\|_{L^4}^{4/3} \ge |\ln \eps|^4) \le c |\ln \eps|^{-2}
$$
it follows that 
 $$
\|w_\eps(0)\|_{L^2}\ds \le c \|v_0\|_{H^2}+ K_{\eps}\|v_0\|_{L^2} |\ln \eps|^4+K_{\eps}\|v_0\|_{H^1}^3
$$
with $K_\eps$ having all moments finite and such that $K_{2^{-k}}$ is almost finite by Borel-Cantelli.We have taken
$|\ln \eps|^4$ instead of $|\ln \eps|^2$ in the estimate above in order to have this latter property.
Recall that $Y_{2^{-k}}$ converges a.s. in $L^\infty$.

We do not have preservation of the $L^2$ norm but: 
$$
\begin{array}{l}
\ds \frac12\frac{d}{dt}\int_{\T^2} |w_\eps(x,t)|^2 e^{-2Y_\eps(x)} dx \\
\\
\ds =2 \lambda \int_{\T^2} Re(v_\eps(x,t) \bar w_\eps(x,t)) Im (v_\eps(x,t) \bar w_\eps(x,t)) e^{-4Y_\eps(x)}dx \\
\\
\ds \le K_{\eps} \|v_\eps(t)\|_{L^\infty}^2 \int_{\T^2} |w_\eps(x,t)|^2 e^{-2Y_\eps(x)} dx\\
\\
\ds \le  K_{\eps} \|v_\eps(t)\|_{H^1}^2\left(1+\ln(1+\|v_\eps(t)\|_{H^2})\right) \int_{\T^2} |w_\eps(x,t)|^2 e^{-2Y_\eps(x)} dx\\
\\
\ds \le  K_{\eps} \left(\|v_0\|_{H^1}^2+\|v_0\|_{L^2}^2\|v_0\|_{H^1}^2\right)\left(1+\ln(1+\|v_\eps(t)\|_{H^2})\right) \int_{\T^2} |w_\eps(x,t)|^2 e^{-2Y_\eps(x)} dx
\end{array}
$$
thanks to the Brezis-Gallouet inequality (\cite{brezis-gallouet}) and to Proposition \ref{p4.1}. Then as above we have:
$$
\begin{array}{ll}
\ds \|\Delta v_\eps(t)\|_{L^2}&\ds \le 2\|w_\eps(t)\|_{L^2}+c K_{\eps}\|v_\eps(t)\|_{L^2} |\ln \eps|^4+\lambda\||v_\eps(t)|^2v_\eps e^{-2Y_\eps}\|_{L^2}\\
\\
&\ds \le 2\|w_\eps(t)\|_{L^2}+ K_{\eps}\|v_\eps(t)\|_{L^2} |\ln \eps|^4+K_\eps \|v_\eps(t)\|_{L^6}^3.
\end{array}
$$
By the embedding $H^1\subset L^6$ and Proposition \ref{p4.1}, we deduce 
$$
\begin{array}{ll}
\| v_\eps(t)\|_{H^2}&\le c \|w_\eps(t)\|_{L^2} + K_{\eps} \left(  \|v_0\|_{L^2} |\ln \eps|^4+\|v_0\|_{H^1}^3+\|v_0\|_{L^2}^3\|v_0\|_{H^1}^3\right).
\end{array}
$$
Again, the almost sure boundedness of the different constant $K_k$ is obtained thanks to Lemma \ref{l2.1bis}, Lemma \ref{l2.1} and Borel-Cantelli.

To lighten the following computation, we use the temporary notations: 
$$
\begin{array}{c}
\tilde w_\eps= \|w_\eps e^{-Y_\eps}\|_{L^2}^2, \quad \alpha_\eps= K_{\eps} \left(\|v_0\|_{H^1}^2+\|v_0\|_{L^2}^2\|v_0\|_{H^1}^2\right),\\
\\
  \beta_\eps= K_{\eps} \left(  \|v_0\|_{L^2} |\ln \eps|^4+\|v_0\|_{H^1}^3+\|v_0\|_{L^2}^3\|v_0\|_{H^1}^3\right).
\end{array}$$
Then we have: 
$$
\begin{array}{ll}
\ds \frac{d}{dt}\tilde w_\eps &\ds  \le \alpha_\eps \left(1+\ln(1+\|v_\eps\|_{H^2})\right) \tilde w_\eps\\
&\ds \le \alpha_\eps \left(1+\ln(1+c\|w_\eps\|_{L^2}+\beta_\eps)\right) \tilde w_\eps\\
&\ds \le  \alpha_\eps \left(1+\ln(1+K_{\eps}\tilde w_\eps+\beta_\eps)\right) \tilde w_\eps.
\end{array}
$$
Hence
$$
\frac{d}{dt}(1+\ln(1+K_{\eps}\tilde w_\eps(t)+\beta_\eps))\le \alpha_\eps (1+\ln(1+K_{\eps}\tilde w_\eps(t)+\beta_\eps)).
$$
By Gronwall's Lemma we deduce:
$$
1+\ln(1+K_{\eps}\tilde w_\eps(t)+\beta_\eps)\le (1+\ln(1+K_{\eps}\tilde w_\eps(0)+\beta_\eps))\exp(\alpha_\eps t)
$$
and taking the exponential
$$
\begin{array}{ll}
\|w_\eps(t)\|_{L^2} &\ds \le K_{\eps}\tilde w_\eps(t) \le c K_{\eps}(1+K_{\eps}\tilde w_\eps(0)+\beta_\eps)^{\exp(\alpha_\eps t)}\\
\\
&\ds \le K_{\eps}(1+K_{\eps} \|w_\eps(0)\|_{L^2}+\beta_\eps)^{\exp(\alpha_\eps t)}\\
\\
&\ds \le K_{\eps}(1+K_{\eps} \|v_0\|_{H^2}+\beta_\eps)^{\exp(\alpha_\eps t)}.
\end{array}
$$
The result follows. 
\hfill$\square$

We see that this $H^2$ bound is not as good as in the linear case. Due to  the double exponential, we do not have moments here. However, this is sufficient to prove existence and uniqueness.  We now state the main result of this section.
\begin{Theorem}\label{t4.1}
Assume that $v_0=u_0e^{Y}\in H^2$ and 
\begin{equation}\label{e4.3}
\|e^{-2Y}\|_{L^\infty}^3\|e^{2Y}\|_{L^\infty} \|v_0\|_{L^2}< 1 \mbox{ if } \lambda=1.
\end{equation}
For any $T\ge 0$, $p\ge 1$, $\gamma\in (1,2)$, when $\eps\to0$, the solution $v_\epsilon$ of \eqref{e3.1} satisfying $v_\eps(0)=v_0$ converges in probability in $C([0,T];H^\gamma))$ to $v$ which is the unique solution to 
\begin{equation}
\label{e4.4}
i\frac{dv}{dt}=\Delta v -2\nabla v\cdot \nabla Y + v :|\nabla Y^2: +\lambda |v|^2v e^{-2Y}
\end{equation}
with paths in $C([0,T];H^\gamma)$ such that $v(0)=v_0$.
\end{Theorem}
{\bf Proof:} Again, pathwise uniqueness is easy. 

Under assumption \eqref{e4.3}, we know that \eqref{e4.2}   holds for $\eps$ small enough. Thus we may use Propositions \ref{p4.1} and \ref{p4.2}.
 
We take $\eps_2>\eps_1>0$, set $r=v_{\eps_1}-v_{\eps_2}$ and write:
$$
\begin{array}{ll}
\ds i\frac{dr}{dt}=&\ds \Delta r-2r\cdot \nabla Y_{\eps_1} + r :|\nabla Y_{\eps_1}|^2:
-2\nabla v_{\eps_2}\cdot \nabla (Y_{\eps_1}-Y_{\eps_2}) + v_{\eps_2} (:|\nabla Y_{\eps_1}|^2: -:|\nabla Y_{\eps_2}|^2:)\\
\\
&\ds+\lambda|v_{\eps_1}|^2re^{-2Y_{\eps_1}}-\lambda (|v_{\eps_2}|^2-|v_{\eps_1}|^2)v_{\eps_2}e^{-2Y_{\eps_1}}
 +\lambda |v_{\eps_2}|^2v_{\eps_2}(e^{-2Y_{\eps_1}}-e^{-2Y_{\eps_2}}).
\end{array}
$$
By the same arguments as in Section \ref{s3} and standard estimates we obtain;
$$
\begin{array}{l}
\ds \frac12\frac{d}{dt}\int_{\T^2}|r(x,t)|^2 e^{-2Y_{\eps_1}(x)}dx \\
\\
\ds \le c\|e^{-2Y_{\eps_1}}\|_{B^{\kappa}_{\infty,2}}\|v_{\eps_2}(t)\|_{H^{1+\kappa}}(\|v_{\eps_1}(t)\|_{H^\kappa}+\|v_{\eps_1}(t)\|_{H^\kappa})\\
\ds\hspace{1.5cm} \times
\left(\|\nabla (Y_{\eps_1}-Y_{\eps_2}) \|_{ B_{\infty,2}^{-\kappa}}+\|:|\nabla Y_{\eps_1}|^2: -:|\nabla Y_{\eps_2}|^2: \|_{ B_{\infty,2}^{-\kappa}}\right)\\
\\
\ds +K_{\eps_1}(\|v_{\eps_1}\|_{L^\infty}+\|v_{\eps_2}\|_{L^\infty})\|v_{\eps_2}\|_{L^\infty}\int_{\T^2}|r(x,t)|^2 e^{-2Y_{\eps_1}(x)}dx\\
\\
\ds + K_{\eps_1} \|v_{\eps_2}\|_{L^\infty}^3(\|v_{\eps_1}\|_{L^2}+\|v_{\eps_2}\|_{L^2})\|e^{-2Y_{\eps_1}}-e^{-2Y_{\eps_2}}\|_{L^2}\\
\\
\ds \le K_{\eps_1,\eps_2}\|v_{\eps_2}(t)\|_{H^{2}}^{\kappa}(\|v_{\eps_1}(t)\|_{H^1}+\|v_{\eps_1}(t)\|_{H^1})^{2-\kappa}\eps_2^{\kappa/2}\\
\\
\ds +K_{\eps_1}(\|v_{\eps_1}\|_{H^1}+\|v_{\eps_2}\|_{H^1})^2(1+\ln(1+\|v_{\eps_1}\|_{H^2})+\ln(1+\|v_{\eps_2}\|_{H^2}))\\
\ds \hspace{1.5cm} \times\int_{\T^2}|r(x,t)|^2 e^{-2Y_{\eps_1}(x)}dx\\
\\
\ds + K_{\eps_1,\eps_2}\|v_{\eps_2}\|_{H^1}^3(1+\ln(1+\|v_{\eps_1}\|_{H^2}))^{3/2}(\|v_{\eps_1}\|_{L^2}+\|v_{\eps_2}\|_{L^2})\eps_2^{\kappa},
\end{array}
$$
by the Brezis-Gallouet inequality.
We have used for $\kappa \in (0,1)$: 
$$
\begin{array}{ll}
\|e^{-2Y_{\eps_1}}-e^{-2Y_{\eps_1}}\|_{L^2}&\le 2 \|Y_{\eps_1}-Y_{\eps_2}\|_{L^2}(\|e^{-2Y_{\eps_1}}\|_{L^\infty}+\|e^{-2Y_{\eps_2}}\|_{L^\infty})\\
\\
&\le c |\eps_2|^{\kappa}\|Y\|_{B_{\infty,\infty}^{\kappa}}(\|e^{-2Y_{\eps_1}}\|_{L^\infty}+\|e^{-2Y_{\eps_2}}\|_{L^\infty})\\
\\
&\le K_{\eps_1,\eps_2} |\eps_2|^{\kappa}.
\end{array}$$
Again, the constants $K_{\eps_1}$, $K_{\eps_1,\eps_2}$ above have all moments bounded independently of $\eps_1,\eps_2$ are almost surely bounded 
in $k$ for $\eps_1=2^{-(k+1)}$, $\eps_2=2^{-k}$.

We deduce from Gronwall's Lemma:
$$
\begin{array}{ll}
\ds \ln \|r\|_{L^2} &\ds\le K_{\eps_1,\eps_2} P(v_0)(1+\ln P(v_0) +\ln|\ln \eps_1|)^{\exp (K_{\eps_1,\eps_2}P(v_0)t)}\\
\\
&\ds+ \kappa \exp(K_{\eps_1,\eps_2}P(v_0)t)(\ln K_{\eps_1,\eps_2}+\ln P(v_0)+ \ln|\ln \eps_1|)-\frac\kappa2 |\ln\eps_2|  ,
\end{array}
$$
where $P(v_0)$ denotes a polynomial in $\|v_0\|_{H^2}$. By interpolation, we have
$$
\ln \|r\|_{H^\gamma}\le c+ (1-\frac\gamma2) \ln \|r\|_{L^2}+\frac\gamma2 \ln \|r\|_{H^2}
$$
so that a similar estimate holds for  $\ln \|r\|_{H^\gamma}$, $\gamma<2$.

The constants $K_{\eps_1,\eps_2}$ are bounded almost surely when $\eps_1=2^{-(k+1)},\; \eps_2=2^{-k}$ so that
$$
\|v_{2^{-(k+1)}}-v_{2^{-k}}\|_{H^\gamma}\le C(v_0) |\ln 2^{-(k+1)}|^{C(v_0)} 2^{-\frac\kappa2(k+1)},
$$
where now $C(v_0)$ is a random constant depending on $\|v_0\|_{H^2}$.  It follows that $v_{2^{-k}}$ is Cauchy in $C([0,T],H^\gamma)$.

Finally, we reproduce the estimate above for $\|v_\eps(t)-v_{2^{-k}}\|_{L^2}$ but bound $\|v_{2^{-k}}\|_{L^\infty}$ by $\|v_{2^{-k}}\|_{H^\gamma}$ instead of using Brezis-Gallouet inequality. 
We obtain:
$$
\begin{array}{ll}
\ds \ln \|v_\eps(t)-v_{2^{-k}}\|_{L^2} &\ds\le \tilde K_{\eps,2^{-k}} P(v_0)(1+\ln P(v_0) +\ln|\ln \eps|)^{\exp (\tilde K_{\eps,2^{-k}}P(v_0)t)}\\
\\
&\ds+ \kappa \exp(\tilde K_{\eps,2^{-k}}P(v_0)t)(\ln \tilde K_{\eps,2^{-k}}+\ln P(v_0)+ \ln|\ln \eps_1|)-\frac\kappa2 |\ln\eps|  ,
\end{array}
$$
where $\tilde K_{\eps,2^{-k}}$ are constants with moments bounded independently on $\eps$ and $k$. Again, a similar bound holds for the $H^\gamma$ norm thanks to an interpolation argument. Letting $k\to \infty$ yields:
 $$
\begin{array}{ll}
\ds \ln \|v_\eps(t)-v\|_{L^2} &\ds\le \tilde K_{\eps} P(v_0)(1+\ln P(v_0) +\ln|\ln \eps|)^{\exp (\tilde K_{\eps}P(v_0)t)}\\
\\
&\ds+ \kappa \exp(\tilde K_{\eps}P(v_0)t)(\ln \tilde K_{\eps}+\ln P(v_0)+ \ln|\ln \eps_1|)-\frac\kappa2 |\ln\eps|  ,
\end{array}
$$
where again $\tilde K_{\eps}$ are constants with bounded moments.

The conclusion follows. \hfill $\square$
\begin{Remark}
Condition \eqref{e4.3} is probably not optimal.   It can easily be weakened to  \newline $\lambda K_{2,\eps}^3K_{1,\eps} \|v_0\|_{L^2}< 2$, but this is probably not optimal either. 
\end{Remark}
\subsection*{Acknowledgement}  

\noindent 
This  work was started during the Fall semester 2015, while the authors were in residence at the Mathematical Sciences Research Institute in Berkeley, California,  supported by the National Science Foundation under Grant No. DMS-1440140. 
A. Debussche  benefits from the support of the French government ``Investissements d'Avenir" program ANR-11-LABX-0020-01. H. Weber is supported by the Royal Society through the University Research Fellowship UF140187.

\bigskip

\noindent Arnaud Debussche, IRMAR, \'Ecole Normale Sup\'erieure de Rennes, UBL, CNRS, Campus de Ker Lann, 37170 Bruz, France. e-mail: \texttt{arnaud.debussche@ens-rennes.fr}.

\bigskip 
\noindent
Hendrik Weber, University of Warwick, Coventry, United Kingdom. \\
e-mail: \texttt{hendrik.weber@warwick.ac.uk}.


\begin{thebibliography}{99}
\bibitem{allez-chouk} R. Allez and K. Chouk, {\it The continuous Anderson Hamiltonian in dimension two}, arXiv:1511:02718v2.

\bibitem{BCD} H. Bahouri, J.-Y. Chemin and R. Danchin, {\it Fourier analysis and nonlinear partial differential equations}, vol. 343 of Grundlehren der Mathematischen Wissenschaften. Springer, Heidelberg, 2011.  
 
 \bibitem{bourgain} J. Bourgain {\it Invariant measures for the 2D-defocusing nonlinear Schr\"odinger equation},  Comm. Math. Phys. 176 (1996), no. 2, pp. 421?445
 
\bibitem{brezis-gallouet} H. Brezis and T. Gallouet
{\it Nonlinear Schr\"odinger evolution equations }, Nonlinear Analysis, Theory, Methods \& Applications, vol. 4, no 4, pp. 677-681.
   
\bibitem{cazenave} T. Cazenave,
{\it Semilinear Schr\"odinger equations}, Courant Lecture Notes in Mathematics, American Mathematical Society,
Courant Institute of Mathematical Sciences, 2003. 
   
 \bibitem{conti1} C. Conti, {\it Solitonization of the Anderson localization}, Phys. Rev. A, {\bf 86}, 2012.
 
 \bibitem{conti2} N. Ghofraniha, S. Gentilini, V. Folli, E. DelRe and C. Conti, {\it Shock waves in disordered media}, Phys. Rev. Letter, {\bf 109}, 2012.
   
\bibitem{HL15} M. Hairer and C. Labb\'e, {\it A simple construction of the continuum parabolic anderson model on $\R^2$}.
Electron. Commun. Probab. 20 (2015), no. 43, pp. 1-11
 

 


  

\end{thebibliography}
\end{document}